\documentclass[12pt]{amsart}
\usepackage{oldlfont}
\usepackage{enumerate}
\usepackage{hyperref}
\usepackage{amssymb}
\usepackage{amsmath}
\usepackage{amsfonts}
\usepackage{color}

\DeclareMathOperator{\lcm}{lcm}

\usepackage[capitalise]{cleveref}
\crefname{corollary}{Corollary}{Corollaries}
\crefname{conjecture}{Conjecture}{Conjectures}
\crefname{proposition}{Proposition}{Propositions}
\crefname{question}{Question}{Questions}
\crefname{lemma}{Lemma}{Lemmas}

\newtheorem{theorem}{Theorem}[section]
\newtheorem{lemma}[theorem]{Lemma}
\newtheorem{notation}[theorem]{Notation}
\newtheorem{proposition}[theorem]{Proposition}

\newtheorem*{hypothesis*}{Main Hypothesis}
\newcommand{\ZZ}{{\rm Z}\kern-3.8pt {\rm Z} \kern2pt}
\newcommand{\RR}{${\rm I\kern-1.6pt {\rm R}}$}
\newcommand{\NN}{{\rm I\kern-1.6pt {\rm N}}}
\newcommand{\HH}{{\rm I\kern-1.6pt {\rm H}}}
\newcommand{\QQ}{${\rm Q}\kern-3.9pt {\rm \vrule height7.1pt
    width.3pt depth-2pt} \kern7.5pt$}
\newcommand{\CC}{${\rm \kern 3.5pt \vrule height 6.5pt
    width.3pt depth-1.3pt \kern-4pt C \kern.8pt}$}
\newcommand{\FF}{${\rm I\kern-1.6pt {\rm F}}$}

\newcommand{\Aut}{{\rm Aut}}

\newcommand{\nl}{\trianglelefteq}
\newcommand{\PSL}{{\rm PSL}}

\newcommand{\PGaL}{{\rm P\Gamma L}}
\newcommand{\GL}{{\rm GL}}
\newcommand{\SL}{{\rm SL}}
\newcommand{\AGL}{{\rm AGL}}
\newcommand{\ASL}{{\rm ASL}}

\newcommand{\Sym}{{\rm Sym}}

\newcommand{\ms}{\medskip}

\renewcommand{\o}{\overline}

\begin{document}

\title{Vertex stabilizers of locally $s$-arc transitive graphs of pushing up type}

\author{John van Bon}
\author{Chris Parker}

\maketitle

\renewcommand\rightmark{{\sc Locally $s$-arc transitive graphs}}
\renewcommand\leftmark{{\sc John van Bon and Chris Parker}}

\pagestyle{headings}

\begin{abstract}
Suppose that $\Delta$  a thick, locally finite and locally $s$-arc transitive $G$-graph with $s \ge 4$.  For a vertex $z$ in $\Delta$, let $G_z$ be the stabilizer of $z$ and $G_z^{[1]}$ be the kernel of the action of $G_z$ on the neighbours of $z$.  We say $\Delta$ is of pushing up type provided there exists a prime $p$ and a $1$-arc $(x,y)$ such that $C_{G_z}(O_p(G_z^{[1]})) \le O_p(G_z^{[1]})$ for $z \in \{x,y\}$ and $O_p(G_x^{[1]}) \le O_p(G_y^{[1]})$. We show that if $\Delta$ is of pushing up type, then  $O_p(G_x^{[1]})$ is elementary abelian and $G_x/G_x^{[1]}\cong X$ with $ \PSL_2(p^a)\le X \le \PGaL_2(p^a)$.
\end{abstract}

\par\vskip .5cm

\section{Introduction}

In this article we consider graphs $\Delta$ that are
connected, undirected and without loops or multiple edges. The vertex set of $\Delta$ is denoted by $V\Delta$ and the edge set is $E\Delta$.
A \emph{$G$-graph} is a graph $\Delta$ together with a subgroup $G \leq \Aut(\Delta)$. An
\emph{$s$-arc} emanating from $x_0\in V\Delta$ is a path $(x_0, x_1, \dots ,x_s)$ with
  $x_{i-1}\neq x_{i+1}$ for $1\leq i \leq s-1$. Denote by  $G_z$ the stabilizer of a vertex $z\in V\Delta$.

A $G$-graph $\Delta$ is
\begin{itemize}
\item[-] {\em thick} if the valency at each vertex is at least 3;
\item[-] {\em locally finite} if for each $z \in V\Delta$, $G_z$ is a finite group;
\item[-] {\em locally $s$-arc transitive} if for every
vertex $z \in V\Delta$, $G_z$ is transitive on the set of $s$-arcs
emanating from $z$.
\end{itemize}

This paper is part of ongoing research aimed at determining all vertex stabilizer amalgams for thick, locally finite and locally $s$-arc transitive $G$-graphs for $s\geq 4$. Throughout this introduction $\Delta$ represents a thick, locally finite $G$-graph.

It is easy to see that if $\Delta$ is a locally $s$-arc transitive $G$-graph with
$s\geq 1$, then
$G$ is transitive on  $E \Delta$ and thus $G$
has at most two orbits  on $V \Delta$. If $s\ge 2$, then, for $x \in V \Delta$, we also know that  $G_x$ acts $2$-transitively on $\Delta(x)=\{v\mid \{x,v\}\in E\Delta\}$.
For  $(x_1,x_2)$   a $1$-arc in   $\Delta$, the triple $(G_{x_1},G_{x_2};G_{x_1}\cap G_{x_2})$ is
called the  \emph {vertex stabilizer amalgam} of $\Delta$ with respect to the $1$-arc $(x_1,x_2)$.
When we study $s$-arc transitive $G$-graphs with $s\ge 1$, it is impossible to determine $\Delta$.  The best we can hope for is a description on the vertex stabilizer amalgam and best of all this will be described up to isomorphism of the amalgam.

For a vertex $z \in V\Delta$, $G_z$  acts on $\Delta(z)$. The kernel of this action is denoted by $G_z^{[1]}$ and $G_z^{\Delta(z)}$ is the permutation group $G_z/G_z^{[1]}$.
A locally finite $G$-graph $\Delta$ is   of
\begin{itemize}
\item[-] {\em local characteristic $p$},
if there exists a prime $p$ such that
$$C_{G_z} (O_p(G_z ^{[1]})) \leq O_p(G_z ^{[1]}) \hbox{, for all } z \in  V\Delta;$$
\item [-]{\em pushing up type} with respect to the $1$-arc $(x ,y)$  and the prime $p$, if
$\Delta$ is of local characteristic $p$   and $$O_p(G_{x} ^{[1]}) \leq O_p(G_{y} ^{[1]}).$$
\end{itemize}

Assume that $(x,y)$ is a $1$-arc and set $G_{x,y}= G_x \cap G_y$.
One consequence of the local characteristic $p$ property is that $O_p(G_x^{[1]})$ and $O_p(G_y^{[1]})$ are non-trivial.  In particular, this means that $G_x$ and $G_y$ are rather large and have potentially complicated structure.  Notice that if $K$ is a subgroup of $G_{x,y}$, and $K$ is normalized by both $G_{x}$ and $G_y$, then $K$ fixes every vertex of $\Delta$ and is   consequently trivial. Hence, if $\Delta$ is of pushing up type with respect to $(x,y)$, then, as $O_p(G_x^{[1]})$ is non-trivial, we learn that $\Delta$ is not of pushing up type with respect to $(y,x)$.  Hence $G$ has two orbits on $V\Delta$.

The generic examples of thick, locally finite  $s$-arc transitive $G$-graphs  $\Delta$ with $s \ge 4$  have vertex stabilizer amalgams which are weak $BN$-pairs \cite{DS}. One of the achievements  in  \cite[Theorem 1]{bonstell} is the proof that, for $s \ge 6$ the generic examples are the only examples. In particular,\cite[Corollary 1]{bonstell} remarks that $s \le 9$ for any such $G$-graph $\Delta$.
In \cite{b3} examples of $G$-graphs which are of pushing up type with $s = 5$ have been constructed via amalgams in $\Sym(p^{2a})$. Thus the vertex stabilizer amalgam of $\Delta$ needs not to be a weak $BN$-pair when $s \le 5$.
To determine the $G$-graphs $\Delta$ which are thick, locally finite and locally $s$-arc transitive with $4\le s \le 5$, as in \cite{bonstell}, we consider three distinct cases:
\begin{enumerate}
\item [-] $\Delta$ is not of local characteristic $p$.
\item [-] $\Delta$ is of local characteristic $p$ but not of pushing up type.
\item [-] $\Delta$ is of pushing up type.

\end{enumerate}

In the first case, the vertex stabilizer amalgams have been completely determined \cite[Theorem 1]{NonLocalChar} and it is shown that $s = 5$.
We expect that the second possibility  yields weak  $BN$-pairs. The configurations appearing in option three are the subject of this article together with its companions  \cite{b1,b2} which are in preparation.

Assume from now on that $\Delta$ is of  pushing up type with respect to $(x,y)$ and the prime $p$.
For  $z \in V\Delta$  set  $Q_z=O_p(G_z^{[1]})$ and
$$L_z = \langle Q_u \mid u \in \Delta (z) \rangle Q_z.$$

The main result of this papers is as follows.

\begin{theorem}\label{thm:main}
Suppose that $s \ge 4$ and $\Delta$ is a thick, locally finite, locally $s$-arc transitive $G$-graph of
pushing up type with respect to the 1-arc $(x,y)$ and the prime $p$.
Then $p$ is odd and the following hold \begin{enumerate}[(a)]\item  $G_x^{\Delta(x)} \cong X$  where $\PSL_2(p^a)\le X \le \PGaL_2(p^a)$ and  $\Delta(x)$ has size $p^a+1$ and is identified with the projective line for $X$;
\item
  $L_x/Q_x \cong \SL_2(p^a)$,  $O^p(L_x) \cong \ASL_2(p^a)^\prime$ and $Q_x$ an elementary abelian $p$-group.\end{enumerate}
\end{theorem}

In \cite{b1,b2} mentioned above, van Bon establishes the isomorphism types of the vertex stabilizer amalgams appearing in the conclusion of \cref{thm:main} and so completes the determination of $G$-graphs of pushing up type with $s \geq 4$.
This then extends \cite[Lemma 7.7]{bonstell} which can be interpreted to say that if $\Delta$ is of pushing up type, then  $s\leq 5$.

 The organization of this paper is as follows.
In Section 2 we derive   general properties of vertex stabilizer amalgams for $\Delta$ of pushing up type. In Section 3, we consider the possibility that $G_x^{\Delta (x)}$ is a projective linear group of degree at least 3,
the main result of the section is \cref{prop:notlnq} which asserts that $F^*(G_x/G_x^{[1]}) \not \cong \PSL_n(p^a)$ with $n \ge 3$. The strategy followed to obtain the conclusion of \cref{prop:notlnq}   uses the results of Section 2 and is similar in flavour to \cite{NonLocalChar} which exploits Zsigmondy primes.
This fails to eliminate the possibility that  $F^*(G_x/G_x^{[1]}) \cong \PSL_3(2)$  however, and so
here we call upon a pushing up result \cite{MS} which allows us to compare non-central chief factors of the vertex stabilizers. This  eventually leads to the elimination of this last case as well.
Finally, in Section 4, we recall that the action of $G_x^{\Delta (x)}$ on $\Delta(x)$ is
$2$-transitive and so with the help of the classification of finite 2-transitive groups and \cref{prop:notlnq} we see that $L_x$ is a rank 1 Lie
type group or is of regular type.  This is precisely the situation handled in
Section 3 of \cite{bonstell}. After application of these results to our case we are left in a situation where we can follow steps 1-9 of the proof of \cite[Lemma 7.7]{bonstell} word for word to obtain the theorem.

Throughout this paper we   assume the following hypothesis:

\begin{hypothesis*} \label{hyp:MH}
The   $G$-graph $\Delta$ is   thick, locally finite, locally $s$-arc transitive with $s \geq 4$ and, in addition, is of
pushing up type with respect to the $1$-arc $(x,y)$ and prime $p$.
\end{hypothesis*}

 The notation used in the paper is standard in the theory of (locally) $s$-arc transitive $G$-graphs and given in Section 2.
Our group theoretic  notation follows \cite{kurz}.

\section{Preliminaries}

In this section we prove some properties of
thick, locally finite, and locally $s$-arc transitive $G$-graphs.
We   assume the Main Hypothesis, though some results
  also hold under weaker assumptions. First we fix the notation used throughout the article.

\par\smallskip\noindent
\begin{notation}\label{not:1stepmore}
Let $d(\cdot , \cdot)$  represent the standard distance function on $\Delta$.
For $u \in V\Delta$, $(u,v)$ a $1$-arc in $\Delta$, $\Theta \subseteq V \Delta$ and $i \ge 1$,
\begin{eqnarray*}
\Delta^{i} (u)&=&\{v \in V\Delta \mid d(u,v)\leq i\};\\
\Delta (u)&=& \Delta^1(u) \backslash \{u\};\\
q_u&=&|\Delta (u)|-1;\\
G_u&=&\langle g\in G \mid u^g=u \rangle; \\
G_\Theta&=& \bigcap_{\theta\in \Theta} G_\theta;\\
 G_u^{[i]}&=&G_{\Delta^i(u)};\\
 G_{\Theta}^{[1]}&=& \bigcap_{\theta \in \Theta} G_\theta^{[1]};\\
Q_u&=&O_p(G_u^{[1]});\\
Z_u&=&\Omega_1(Z(Q_u));\\
C_u&=&\langle G_v^{[2]}\mid v \in \Delta(u)\rangle Q_u;\\
L_u &=& \langle Q_v \mid v\in \Delta (u) \rangle Q_u;\\
L_{u,v}&=& G_{u,v} \cap L_u;
\end{eqnarray*}
and $(w,x,y,z)$ is a fixed $3$-arc where $(x,y)$ is an arc for which $\Delta$ is of pushing up type.
\end{notation}
\par\smallskip\noindent

Notice that we do not know that $L_{u,v} = L_{v,u}$ and so the order of the vertices  on the arc is important for the definition of $L_{u,v}$. We recall the fundamental properties of the $G$-graph $\Delta$:
\begin{eqnarray*}C_{G_u}(Q_u)&\le& Q_u\text{ for } u \in \{x,y\}; \text{ and }\\Q_x&\le& Q_y.
\end{eqnarray*}

An essential tool when studying such $G$-graphs is given by the following lemma:

\begin{lemma}\label{lem:basic1}
Assume that $R\le G_{x,y}$, $N_{G_x}(R)$ is transitive on $\Delta(x)$ and $N_{G_y}(R)$ is transitive on $\Delta(y)$. Then $R=1$.
\end{lemma}

\begin{proof} See \cite[10.3.3]{kurz}.  \end{proof}

Because  $\Delta$ has pushing up type,  $Q_x \le Q_y$. In fact, $Q_x< Q_y$ as otherwise $Q_x=Q_y=1$ by \cref{lem:basic1} and then $\Delta$ is not of local characteristic $p$. One consequence of this observation is that $G$ has two orbits on $\Delta$.

We   present  some elementary properties of various of the subgroups defined above.
\begin{lemma}\label{lem:P1}
Let $N_x=O^p(L_{x})$,
and $\widetilde{G}_{x}=G_{x}/Q_{x}$. Then
\begin{enumerate}[(i)]
\item $[ G_{x}^{[1]},L_{x}]\leq Q_{x}$;
\item $N_x \not\leq G_{x}^{[1]}$;
\item  $L_{x}= N_xQ_{y}$ and $G_x=N_xG_{x,y}$;
\item  either  $\widetilde{N_x}$ is quasisimple or an $r$-group, $r$ a prime  with $r\ne p$;
\item $\widetilde{G_x^{[1]}\cap L_x}=\Phi(\widetilde{L}_{x})=Z(\widetilde{L}_{x})$;
\item if $\widetilde{N_x}$ is an $r$-group,
then $\widetilde{G_x^{[1]}\cap L_x}=\widetilde{N_x}^\prime = \Phi(\widetilde{N_x})$,
and $G_{x}$ acts transitively on the non-trivial elements of $N_x^{\Delta (x)}$.
\end{enumerate}
\end{lemma}

\begin{proof}
 The first statement follows from \cite[Lemma 5.1 (a)]{bonstell} while the remainder of the statements can be found in \cite[Lemma 5.2]{bonstell}.
\end{proof}

\begin{lemma}\label{lem:P2}
The following hold:
\begin{enumerate}[(i)]
\item $G_{x,y}^{[1]}=G_x^{[2]}$, $G_{y}^{[2]}=G_y^{[3]}$ and $Q_x=O_p(G_x^{[2]})$;
\item $G_{x} ^{[1]} \cap G_{z} = G_{x}^{[2]}$.
\end{enumerate}
In particular, $G_{x} ^{[1]}$ induces a semi-regular group
on $\Delta (y) \backslash \{{x}\}$, $|G_{x}^{[1]}G_y^{[1]}/G_{y}^{[1]}|$ divides $q_y$ and
$|G_{x,y,z}| = |G_{x,y,z} ^{\Delta (x)}||G_{x} ^{[2]}|.$
\end{lemma}

\begin{proof}
$(i)$.
By \cref{lem:P1} (i) and (iii), $[L_x,G_x^{[1]}] \leq Q_x$  and $L_x$ acts  transitively on $\Delta (x)$. Hence
$$[L_x,G_{x,y}^{[1]}]\leq [L_x,G_x^{[1]}]  \leq Q_x \leq G_{x,y}^{[1]}.$$ Thus $L_x$ normalizes $G_{x,y}^{[1]}$. The
transitivity of $L_x$ on $\Delta (x)$ now yields $G_{x,y}^{[1]}=G_x^{[2]}$.
In particular $G_y^{[2]} \leq G_{u,y}^{[1]}=G_u^{[2]}$ for any $u\in \Delta (y)$.
It follows that $G_{y}^{[2]}=G_y^{[3]}$. Since $Q_x \leq Q_v \leq G_v^{[1]}$ for all $v \in \Delta (x)$, we also have
$Q_x \nl G_x^{[2]}$. Hence $Q_x=O_p(G_x^{[2]})$.
\ms
\par\noindent
$(ii)$.
By \cref{lem:P1} (i)
$$[L_x , G_x^{[1]}\cap G_z]\le [L_x , G_x^{[1]}] \leq Q_x \leq  G_x^{[1]}\cap Q_y\leq G_x^{[1]}\cap G_z.$$
Thus $L_x$ normalizes $G_x^{[1]}\cap G_z$.
Pick  $g \in L_x$  with
$y^g \not = y$ and put $y^\prime = y^g$ and $z^\prime = z^g$.
Then $G_x^{[1]}\cap G_z= G_x^{[1]}\cap G_{z^\prime}$.
Since $s\ge 4$,  $G_{{z^\prime},{y^\prime},x,y}$ acts transitively on $\Delta (y) \backslash \{x\}$ and it also normalizes
$G_x^{[1]}\cap G_{z^\prime}= G_x^{[1]}\cap G_z$. Therefore,
  $ G_x^{[1]} \cap G_{z} \leq G_{y} ^{[1]}$.
Hence, from (i), $$G_x^{[2]} \leq G_x^{[1]}\cap G_z \leq  G_{x,y} ^{[1]}=G_x ^{[2]}.$$
Now $(ii)$ follows.
\ms

Finally, as
$(G_{x} ^{[1]} \cap G_{z})^{\Delta(y)} = 1$ by (ii),
$G_{x} ^{[1]}$ induces a group which acts semi-regularly
on $\Delta ({y}) \backslash \{x\}$. Therefore, we also have $|G_{x} ^{[1]}G_y^{[1]}/G_y^{[1]}|$ divides $\Delta(y)-1=q_y$. To see that $|G_{x,y,z}| = |G_{x,y,z} ^{\Delta (x)}||G_{x} ^{[2]}|,$ just observe that (ii) gives  $$G_x^{[2]} \leq G_{x}^{[1]} \cap G_{x,y,z} \leq G_{x}^{[1]}\cap G_z = G_{x}^{[2]}.$$
\end{proof}

\begin{lemma}\label{lem:P3}
The following hold:
\begin{enumerate}[(i)]
\item $G_y^{[2]}$ and $G_x^{[1]} \cap G_z ^{[1]}$ are   $p$-groups;
\item $G_x^{[1]} \cap G_z ^{[1]} =G_x^{[2]} \cap G_z ^{[2]}= G_{x,y,z}^{[1]}=Q_x \cap Q_z$;
\item $Q_yG_x^{[1]} \cap G_z^{[2]}=Q_z$;
\item $Q_yG_{x}^{[2]} \cap Q_yG_{z}^{[2]} = Q_y$;
\item either $G_x^{[2]}/Q_x$ is abelian or $G_x^{[1]}=G_x^{[2]}$;

\item   ${\mathrm lcm}(q_x,q_y)$ divides $|G_{w,x,y,z}|$  and $\pi(G_{x,y})=\pi (G_{w,x,y,z})$.
\end{enumerate}
\end{lemma}

\begin{proof}
$(i)$.
Let $R \in {\mathrm Syl}_r (G_x^{[1]})$, with $r \neq p$.
 Then \cref{lem:P1} (i) yields $[R ,L_x] \leq [G_x^{[1]},L_x]  \le Q_x$. Hence $L_x$ acts on $Q_xR$ and the Frattini argument gives $L_x= Q_xN_{L_x}(R)$. In particular, $N_{L_x}(R)$ is transitive on $\Delta(x)$.
Now using \cref{lem:P1} (i) we get $[R ,N_{L_x}(R)] \leq [G_x^{[1]},L_x]\cap R\le  Q_x \cap R=1$.
 So for all $R_0 \le R$, $C_{L_x}(R_0)\ge C_{L_x}(R)=N_{L_x}(R)$ is transitive on $\Delta (x)$.

Let $T \in {\mathrm Syl}_t (G_y^{[2]})$, with $t \neq p$. By the Frattini argument
$G_y= N_{G_y}(T)G_y^{[2]}$. Hence $N_{G_y}(T)$ is transitive on $\Delta (y)$.
Since $T \leq G_x^{[1]}$ we have
$C_{L_x}(T)$ is transitive on $\Delta (x)$. Thus
$T$ is normalized by both $N_{G_y}(T)$ and $C_{L_x}(T)$ and so $T=1$ by \cref{lem:basic1}. Thus $G_y^{[2]}$ is a $p$-group as claimed.

Suppose that $T\in {\mathrm Syl}_t(G_x^{[1]} \cap G_z^{[1]})$, with $t \neq p$.
Since $s \geq 4$, there exists a $c \in C_{L_x}(T)$ with $y^c \neq y$ and $z^c \neq z$.
Let $\gamma= (z^c,y^c,x, y)$ and $N=G_x^{[1]} \cap G_{z^c}^{[1]}$.
Now $T=T^c \in {\mathrm Syl}_t(N)$  and $N\trianglelefteq G_{\gamma}$. Then, by \cref{lem:P2} (ii), $N \leq G_x^{[2]}\leq G_y^{[1]}$.
Since $s \geq 4$, $G_{\gamma}$ is transitive on $\Delta (y) \backslash \{x\}$.
By the Frattini argument
$G_{\gamma}= N_{G_{\gamma}}(T)N$. Thus $N_{G_{\gamma}}(T)$
is transitive on $\Delta (y) \backslash \{x\}$ since $N \leq G_y^{[1]}$. As $T\le G_{z}^{[1]}$,
$T \leq G_u^{[1]}$, for all $u \in \Delta (y)$. But then
$T \leq G_y^{[2]}$ which is a $p$-group and thus $T=1$. This completes the proof of (i).
\ms

\par\noindent
$(ii)$. From \cref{lem:P2} (ii) we have  $G_x^{[1]} \cap G_z =G_x^{[2]}$ and $G_x \cap G_z^{[1]}=G_z^{[2]}$ and thus
$$G_x^{[2]} \cap G_z^{[2]} \leq G_{x,y,z}^{[1]} \leq G_x^{[1]}\cap G_z^{[1]} \leq G_x^{[2]} \cap G_z^{[2]}$$ which gives
 $G_x^{[1]}\cap G_z^{[1]}= G_{x,y,z}^{[1]}=G_x^{[2]} \cap G_z^{[2]}$.

Since, by (i), $G_x^{[1]} \cap G_z^{[1]}=G_x^{[2]} \cap G_z^{[2]}$ is a $p$-group which is
normalized by $G_x^{[2]}$ and $G_z^{[2]}$, \cref{lem:P2} (i) implies
$G_x^{[1]} \cap G_z^{[1]} \leq O_p (G_x^{[2]})=Q_x$ and
$G_x^{[1]} \cap G_z^{[1]} \leq O_p (G_z^{[2]})= Q_z$.
Hence $$Q_x \cap Q_z \leq G_x^{[2]} \cap G_z^{[2]} \leq Q_x \cap Q_z$$ and this completes the proof of (ii).
\ms

\par\noindent $(iii)$.  Suppose that $R\le Q_yG_x^{[1]} \cap G_z^{[2]}$ has $p'$-order. Since $G_x^{[1]}$ is normalized by $Q_y$ we have
$R \leq G_x^{[1]}$. Hence
$R \leq G_x^{[1]} \cap G_z^{[2]}$ which by (i) is a $p$-group. Thus $R=1$ and
$Q_yG_x^{[1]} \cap G_z^{[2]}$ is a $p$-group. Since $Q_z \leq Q_yG_x^{[1]} \cap G_z^{[2]} \nl G_z^{[2]}$
and $O_p(G_z^{[2]})=Q_z$ we have
$Q_yG_x^{[1]} \cap G_z^{[2]}=Q_z$.
\ms
\par\noindent
$(iv)$. Using (iii) and the modular law we have $$
Q_yG_x^{[2]} \cap Q_yG_z^{[2]}=Q_y(Q_yG_x^{[2]}\cap G_{z}^{[2]}) \leq Q_y(Q_yG_x^{[1]} \cap G_z^{[2]})=Q_yQ_z=Q_y.$$

\ms
\par\noindent
$(v)$.
 Suppose $G_x^{[1]} \neq G_x^{[2]}$.
Let $g \in G_x^{[1]} \setminus G_x^{[2]}$. Then, by \cref{lem:P2} (ii), $z^g \neq z$. Then $G_x^{[1]}G_z^{[2]}=(G_x^{[1]}G_z^{[2]})^g=G_x^{[1]}G_{z^g}^{[2]}$. Hence, using (ii) and the modular law, we have
\begin{eqnarray*}
[G_z^{[2]},G_z^{[2]}]&\leq &[G_x^{[1]}G_z^{[2]},G_x^{[1]}G_z^{[2]}] \cap G_z^{[2]}\leq [G_x^{[1]}G_z^{[2]},G_x^{[1]}G_{z^g}^{[2]}]\cap G_z^{[2]}
\\&\leq& G_{x}^{[1]}[G_z^{[2]},G_{z^g}^{[2]}]\cap G_z^{[2]}\leq G_{x}^{[1]}(G_z^{[2]} \cap G_{z^g}^{[2]})\cap G_z^{[2]} \\&\leq& G_{x}^{[1]} Q_z\cap G_z^{[2]}=  Q_z.
\end{eqnarray*}
It follows that $G_z^{[2]}/Q_z$ is abelian.

\ms
\par\noindent
$(vi)$. Since $s \geq 4$, $|G_{w,x,y,z}|$ is divisible by $\lcm(q_xq_y)$.
Let $t \in \pi (G_{x,y})$ and $h \in G_{x,y}$ have order $t$. If $t$ divides $q_xq_y$ then $t$ divides $|G_{w,x,y,z}|$.
Suppose $t$ does not divide $q_xq_y$. Then $\langle h\rangle$   fixes a vertex in both
$ u^\prime \in\Delta (x) \backslash \{y\}$ and $ z^\prime \in\Delta (y) \backslash \{x\}$. Thus $t$ divides $|G_{u^\prime,x,y,z^\prime}|=|G_{w,x,y,z}|$ since $s \geq 4$.
It follows that $\pi (G_{x,y}) \subseteq \pi (G_{w,x,y,z})$. The reverse inclusion is
immediate since $G_{w,x,y,z}$ is a subgroup of $G_{x,y}$.
\end{proof}

\begin{lemma}\label{lem:Cy}
If $G_x^{[2]}/Q_x$ is abelian of exponent $t$, then $C_y/Q_y$ is abelian of exponent $t$ and order at least $t^2$.
\end{lemma}

\begin{proof} For  $u,v \in \Delta(y)$ $[G_u^{[2]},G_v^{[2]}] \leq  Q_u\cap Q_v\le Q_y$ by \cref{lem:P3} (ii) and the fact that $G_x^{[2]}/Q_x$ is abelian. Hence $C_y/Q_y$ is abelian and, as $G_x^{[2]}Q_y\cap G_u^{[2]}Q_y= Q_y$ by \cref{lem:P3} (iv), the claim follows.
\end{proof}

The proof of the following lemma is based on an argument that can be found in \cite[Lemma 4.8]{bonstell}.

\begin{lemma}\label{lem:P4}
Let $t$ be a prime dividing $q_x$.
Then $G_x^{[2]}$ is a $t^\prime$-group if and only if $G_x^{[1]}$ is a $t^\prime$-group.
\end{lemma}

\begin{proof}It suffice to prove that $G_{x}^{[1]}$ is a $t'$-group whenever $G_{x}^{[2]}$ is a $t'$-group.  Assume that $G_{x}^{[2]}$ is a $t'$-group.
Let $(x,y,z,a)$ be a $3$-arc and let $T \in {\mathrm Syl}_t(G_{x,y,z,a})$. Then, as $s\ge 4$ and $t$ divides $q_x$, $T$  acts non-trivially on $\Delta(x)$. Hence, again as $s \ge 4$, $T$ fixes no $4$-arcs starting with a vertex in $a^G$. If $N_{G_y}(T) \not \le G_z$, then there exists $g\in N_{G_y}(T)$ such that $T$ fixes the $4$-arc $(a,z,y,z^g,a^g)$, which is a contradiction. Therefore $N_{G_y}(T)\le G_z$. Let $X \in {\mathrm Syl}_t(G_{x}^{[1]})$ be normalized by $T$.  Then $N_{X}(T)\le G_z$. Hence $$ N_{X}(T) \le X \cap G_z\le G_{x}^{[1]} \cap G_z = G_x^{[2]}$$  by \cref{lem:P2} (ii) and so $N_X(T)=1$ as $G_x^{[2]}$ is a $t'$-group. It follows that $ C_{X}(T) \le N_X(T)=1$ and so we conclude that $X=1$. Hence $G_x^{[1]}$ is a $t'$-group and this   establishes the claim.
\end{proof}

The next lemma, which is fundamental for the approach to our proof of \cref{thm:main}, is the origin of the name pushing up type amalgam.

\begin{lemma}\label{lem:no char sbgrp} Suppose that $K$ is a non-trivial characteristic subgroup of $Q_y$. Then $N_{G_{x}}(K)= G_{x,y}$. In particular,
no non-trivial characteristic subgroup of $Q_y$ is normalized by $L_x$.
\end{lemma}

\begin{proof} Suppose $K$ is a non-trivial characteristic subgroup of $Q_y$. Then $K$ is normalized by $G_y$ and so also by $G_{x,y}$. Since $G_x$ acts $2$-transitively on $\Delta(x)$, $G_{x,y}$ is a maximal subgroup of $G_x$. Hence if $N_{G_x}(K) > G_{x,y}$, then $K$ is normalized by $G_x$. But then $K$ is normalized by $G_x$ and $G_y$ and so $K=1$ by \cref{lem:basic1}.
 \end{proof}

\section{$L_x^{\Delta(x)}$ is not a projective linear group of degree at least $3$}

In this section we intend to demonstrate

\begin{proposition}\label{prop:notlnq}
Suppose that   $s\ge 4$ and $\Delta$ is a thick, locally finite, locally $s$-arc transitive $G$-graph of
pushing up type with respect to $(x,y)$ and the prime  $p$. Then, for $n \ge 3$ and $a$ a natural number,
 $F^*(G_x^{\Delta (x)}) \not \cong \PSL_n(p^a)$ acting on projective points.
 \end{proposition}

Throughout this section, we assume that the Main Hypothesis holds and
$$L_x ^{\Delta(x)} \cong \PSL_n(q)$$
with $n \geq 3$, $q=p^a$ and $\Delta(x)$ corresponding to the points in projective $(n-1)$-space. We continue with the notation established in \cref{not:1stepmore}.

 Before we start on the proof, we record the following facts about projective linear groups.

\begin{lemma}\label{lem:P5}
Assume that $n \ge 2$, $p$ is a prime, $q=p^k$ and $\PSL_n(q) \nl H \leq \PGaL_n(q)$ acting on the
projective space $PV$. Let    $u,v \in PV$ be   distinct points.
The following hold:
\begin{enumerate}[(i)]
\item $|H/H_v|-1= q\left (\frac{q^{n-1}-1}{q-1}\right)$.
\item There exists a unique $E\nl H_{v}$ such that
$O_p(E)=O_p(H_v)$, $E/O_p(E) \cong \SL_{n-1}(q)$
and $|H_v/ E|$ divides $(q-1)k$.
Moreover, $O_p(E)$ is a natural module for $E/O_p(E)$.
\item Either $n=3$ and $|H_{u,v}/O_p( H_{u,v})|$
divides $(q-1)^2k$ or there exists  $O_p(H_{u,v}) \nl F\nl H_{u,v}$ such that
$F/O_p(F) \cong \SL_{n-2}(q)$ and
$|H_{u,v}/F|$ divides  $(q-1)^2k$.
\item
Let $N$ and $E$ be subgroups of $H_v$ with $O_p(H_v)=O_p(N)=O_p(E)$ and $E/O_p(E) \cong \SL_{n-1}(q)$.
Suppose $E$ normalizes $N$.
Then one of the following hold:
\begin{enumerate}
\item $E \leq N$;
\item $N/O_p(H_v)$ cyclic, $[N,E] \leq O_p(H_v)$ and $|N/O_p(H_v)|$ divides $(q-1)$;
\item $E^\prime \leq N < E$, $n=3$ and either $N/O_p(H_v) \cong {\mathrm C}_3$ and $q=2$, or $N/O_p(H_v) \cong {\mathrm  Q}_8$ and $q=3$.
\end{enumerate}
\end{enumerate}
\end{lemma}

\begin{proof}
$(i)-(iii)$. Straight forward.
\par\noindent
$(iv)$. For a subgroup $X \leq H_v$ we will denote with
$\o X$ its image in $H_v/(O_p(H_v))$. Since $E$ normalizes $N$ we have
$[\o N, \o E] \nl \o E$.
Hence either $\o E \leq [\o N,\o E]$,  $[\o N,\o E] \leq Z(\o E)$,
or $n=3=q$ and $[\o N, \o E] = [\o E, \o E] \cong Q_8$
or $n=3$, $q=2$ and $[\o N, \o E]= [\o E, \o E] \cong C_3$.
In the second case the Three Subgroup Lemma gives
$[\o N ,[\o E, \o E]]=1$. It follows that
$\o N \leq C_{\o H_v}( [\o E, \o E])$.
Thus  $[\o N, \o E]=1$ too and $\o N$ is a cyclic group whose order divides $q-1$.
In the third case $\o E \leq \o N$ or $\o N \cong {\mathrm Q}_8 $.
In the fourth case $\o E \leq \o N$ or $\o N \cong {\mathrm C}_3$.
\end{proof}

We begin the proof of \cref{prop:notlnq} with a lemma which restricts the structure of $G_x^{[2]}/Q_x$.

\begin{lemma}\label{lem:A2} One of the following holds: \begin{enumerate}
\item[(i)] $G_x^{[2]}/Q_x \cong {\mathrm C_t}$ with $t$ dividing $q-1$; or
\item[(ii)] $ G_x ^{\Delta(x)}  \cong \PSL_3(2)$ and $G_x^{[2]}/Q_x \cong {\mathrm C}_3$.
\end{enumerate}
\end{lemma}

\begin{proof}
Let $P=G_{x,y} $ and put $\overline {P} = P/G_{x}^{[1]}Q_y$. Then, as $G_x^{[1]}Q_y= G_x^{[1]}O_p(P)$, \cref{lem:P5} implies
$\SL_{n-1}(q) \nl \overline P \leq  \Gamma {\mathrm L}_{n-1}(q)$. Let $\overline E \leq \overline P$ with $\overline E\cong \SL_{n-1}(q)$.
Recall that $C_y \le G_y^{[1]}$ and $C_y$ is normal in $G_y$. Hence $\overline C_y$ is a normal subgroup of $\overline P$.

Let $  u \in \Delta(y)\setminus\{x\}$ with $z\ne u$. Then, by \cref{lem:P3} (iv), $[G_{z}^{[2]},G_{u}^{[2]}]\leq Q_y$. It follows that $Q_yG_z^{[2]}$ is normal in $C_y$. Since $\overline {G_z^{[2]}} $ is normal in $\overline {C_y}$, $\overline {G_z^{[2]}}$ is subnormal in $\overline P$.
If $\overline{G_z^{[2]}}$ does not centralize $\overline E$, then $\overline E'\le \overline{G_z^{[2]}} $ by \cref{lem:P5} (iii) and \cref{lem:Cy}.
Since this is true for all $z \in \Delta(y)\setminus\{x\}$ and  $\overline {[G_{z}^{[2]},G_{u}^{[2]}]}=1$, this is impossible unless $(n,q)=(3,2)$ and $\overline{G_z^{[2]}} $ is cyclic of order $3$.  Hence, in the general case, $ \overline{G_z^{[2]}}$ centralizes $\overline {E} $ and we conclude that $\overline{G_z^{[2]}}$  is cyclic of order $t$ dividing $q-1$. We have demonstrated
$$\overline{G_z^{[2]}}\cong \begin{cases}{\mathrm C}_t& t \text{ divides } q-1\cr \mathrm{C}_3&(n,q)= (3,2).\end{cases}$$
Finally, \cref{lem:P3} (iii) yields
$$G_z^{[2]}/Q_z=G_z^{[2]} /(G_z^{[2]}\cap G_x^{[1]}Q_y)\cong    G_z^{[2]}G_{x}^{[1]}Q_y/G_{x}^{[1]}Q_y=\overline{G_z^{[2]}}$$ and this completes the proof.\end{proof}

\begin{lemma}\label{lem:A2.5}
If
$F^*(G_x ^{\Delta(x)}) \cong \PSL_n(q)$  with $n \geq 3$, then $F^*(G_x ^{\Delta(x)}) \cong \PSL_3(2)$ and $G_x^{[2]}/Q_x \cong {\mathrm C}_3$.
\end{lemma}

\begin{proof} Suppose that $F^*(G_x ^{\Delta(x)}) \not  \cong \PSL_3(2)$ with $G_x^{[2]}/Q_x \cong {\mathrm C}_3$.
We first prove that

\ms
\noindent{${\mathbf 1}^\circ$.}  $q^{n-1}  = 2^6$ or $q^{n-1}-1 = p^2-1$ and $p$ is a Mersenne prime.

\ms
\par\smallskip
Suppose that $t$ is a Zsigmondy prime for $((n-1)k,p)$.
By \cite[3.8]{NonLocalChar} $t$ does not divide $(n-1)k$, and
since $n-1 \geq 2$, $t$ does not divide $q-1$.
 Therefore, as $t$ divides $q_x=q(q^{n-1}-1)/(q-1)$, combining \cref{lem:P4,lem:A2} yields $t$ does not divide $|G_x^{[1]}|$. Since $t$ divides $|G_{x,y}|$, $t$ divides $|G_{w,x,y}|$ by \cref{lem:P3} (vi). Therefore  $t$ divides $|{G_{w,x,y}}^{\Delta(x)}|$. However,  \cref{lem:P5} (iii) then implies
 $t$ divides $| \PSL_{n-2}(q)|$, contrary to $t$ being
Zsigmondy prime.  Therefore, Zsigmondy's Theorem (see \cite[3.8]{NonLocalChar})
 implies that
 $q^{n-1}  = 2^6$ or $q^{n-1}-1 = p^2-1$ and $p= 2^r-1$ is a Mersenne prime.

\ms
\noindent{${\mathbf 2}^\circ$.}  We have $q^{n-1}\ne 2^6$.
 \ms

 Assume that $q^{n-1}= 2^6$.  Then, as $n \ge 3$, we have one of the following cases $F^*(G_x^{\Delta(x)})\cong \PSL_3(8)$ with $q_x= 8(8+1)$, $\PSL_4(4)$ with $q_x= 4(4^2+4+1)$ or $\PSL_7(2)$ with $q_x= 2(2^6-1)$. Since $s\ge 4$, $G_{x,y}$ has order divisible $q_x^2$. If the first case  $|G_x^{[1]}|$  is coprime  to $3$ by \cref{lem:P4,lem:A2}. Therefore $|G_{x,y}|_3 =|G_{x,y}^{\Delta(x)}|_3=3^3< 3^4=(q_x^2)_3$, which is a contradiction. Similarly, the second case is impossible as $7^2$ does not divide  $|G_{x,y}^{\Delta(x)}|$. Therefore, $F^*(G_x^{\Delta(x)})\cong \PSL_7(2)$ and $G_x^{[2]}=Q_x$ is a $2$-group by \cref{lem:A2}. Set  $H=G_{w,x,y,z}$. Then by \cref{lem:P3} (v)
$\pi (H)=\pi (G_{x,y})\supseteq \{2,3,5,7,31\}$ and by \cref{lem:P2} (ii) $H\cap G_{x}^{[1]} \le G_z\cap G_{x}^{[1]}$ is a $2$-group. Hence $H^{\Delta(x)}\le G_{w,x,y}^{\Delta (x)} \cong 2^{10}: \SL_5(2)$ and  $\pi (H^{\Delta(x)}) \supseteq \{3,5,7,31\}$. By \cite{Kantor} the maximal over-groups of a Singer cycle in $\SL_5(2)$ have order $5.31=155$ and so we conclude that $G_{w,x,y}^{\Delta (x)}= H^{\Delta (x)}O_2(G_{w,x,y}^{\Delta (x)})$. In particular, we see that $H$ has a quotient isomorphic to $\SL_5(2)$ and has order $2^{\ell}.3^2.5.7.31$ for some $\ell$. Since $s \ge 4$, $H$ operates transitively on $\Delta(z) \setminus \{y\}$. Hence for $a\in  \Delta(z) \setminus \{y\}$, $|H:H_a|= q_x=126$ and so $|H_a|= 2^{\ell-1}.5.31$. Therefore $$ 2^{10}.5.31\ge |H_a^{\Delta(x)}O_2(G_{w,x,y}^{\Delta (x)})/O_2(G_{w,x,y}^{\Delta (x)})|\ge 2^9.5.31>5 .31$$ and contains a Singer cycle of $\SL_5(2)$, which is a contradiction.

\ms

Because of $({\mathbf 1}^\circ)$ and $({\mathbf 2}^\circ)$, it remains to exclude the possibility that $q^{n-1}-1 = p^2-1$ with $p=2^r-1$ a Mersenne prime. Thus, $n=3$, $q_x = p(p+1) = 2^rp$ and $p-1=2(2^{r-1}-1)$ is not divisible by $4$.

 By  \cref{lem:P2} (ii), $G_{w,x,y,z} \cap G_x^{[1]}=G_{x}^{[1]} \cap G_z=G_x^{[2]}$ and by \cref{lem:A2} $|G_x^{[2]}|$ divides $(p-1)|Q_x|$.
By \cref{lem:P5} $|G_{w,x,y,z}^{\Delta (x)}|$ divides $(p-1)^2p^2$. Hence $|G_{w,x,y,z}|$ divides $(p-1)^3p^2|Q_x|$. So, as $4$ does not divide $(p-1)$, and $G_{w,x,y,z}$ acts transitively on $\Delta(z)\setminus\{y\}$, $q_z=q_x=2^rp$ is not divisible by $16$. Hence $r \in\{2,3\}$. Furthermore, if $r=3$, then $|G_x^{[2]}|$ is even. We signal

\ms

\noindent{${\mathbf 4}^\circ$.} $q=p=2^{r}-1\in \{3,7\}$  and $|G_x^{[2]}|$ is even if $q=7$.

\ms

Suppose  that $G_x^{[2]}$ has odd order. Then $r=2$, $q=p=3$, $G_x=G_x^{[1]}L_x$, $L_x/Q_x \cong \PSL_3(3)$ and $q_x=12$.
By \cref{lem:P4}, $G_x^{[1]}$  also has odd order.
Let $S \in {\mathrm Syl}_2(G_{w,x,y,z})$. Then $S \cap G_x^{[1]}=1$ and so $S\cong S^{\Delta(x)} \leq G_{w,x,y}^{\Delta(x)}$.
Hence $S$ is elementary abelian and $|S|\leq 4$.
Since $q_x=12$ and $q_x$ divides $|G_{w,x,y,z}|$ by \cref{lem:P3} (v), we conclude that
$|S|=4$. Since $s \ge 4$, and $q_z=q_x=12$, we now know that $|G_{w,x,y,z,a}|$  is odd  for all $a \in \Delta(z)\setminus\{y\}$.

Let $T \in {\mathrm Syl}_2(G_{y,z})$ with $S \leq T$.
Then $T \cap G_z^{[1]}=1$ and $T \cong {\mathrm SDih}(16)$ is semidihedral.
Since $S$ is elementary abelian of order 4, we have $Z(T) \leq S$. In particular, $Z(T)G_{z}^{[1]}/G_{z}^{[1]}$ inverts $Q_yG_{z}^{[1]}/G_{z}^{[1]}$.
Observe that $Q_y$ has $4$ orbits of length
$3$ on $\Delta (z) \setminus \{y\}$ each of which is fixed by
$Z(T)$. Hence $Z(T)$ fixes a vertex $a \in \Delta(z) \backslash \{y\}$, contrary to
$|G_{w,x,y,z,a}|$ being odd. This contradiction shows that  $|G_x^{[2]}|$ is divisible by 2.

\ms

Assume that $|G_x^{[2]}|$ is even. Then $q\in \{3,7\}$ by $({\mathbf 4}^\circ)$. \cref{lem:A2} states that $G_x^{[2]}/Q_x$ is cyclic and so the Sylow 2-subgroups of $G_x^{[2]}/Q_x$ have order $2$. By \cref{lem:Cy}, $C_y /Q_y$ is abelian and the Sylow $2$-subgroups of $C_y$  are elementary abelian and  have order at least $4$. As $ C_y ^{\Delta{x}}$ is normal in $G_{x,y}^{\Delta(x)}$ and $G_{x,y}^{\Delta(x)}/Q_y^{\Delta(x)}$ is isomorphic to a subgroup of $\GL_2(q)$ containing $\SL_2(q)$, \cref{lem:P5}(iv) yields that
$C_y^{\Delta (x)}/Q_y^{\Delta(x)}$ is cyclic.  \cref{lem:P2} (ii) gives $G_x^{[1]} \cap C_y\leq G_x^{[1]} \cap G_y^{[1]}=G_x^{[2]}$ and thus
$C_y^{\Delta (x)}/Q_y^{\Delta(x)}\cong C_y/G_x^{[2]}Q_y$.
We deduce that $C_y/Q_y$
contains exactly three involutions.  If $q_y >2$, then there exists $u, v \in \Delta(y)$ such that $G_u^{[2]}Q_y/Q_y \cap G_{v}^{[2]}Q_y/Q_y$ has an involution and this contradicts \cref{lem:P3} (iv).  Hence $q_y=2$ and $G_y$ acts transitively on the three involutions in $C_y/Q_y$. Let $S_y \in {\mathrm Syl}_2(C_y)$. Then $S_yQ_y$ is normalized by $G_y$ and thus $S_yQ_y$ has a unique $G_y$-conjugacy class of involutions.

Since $q_y=2$, we have $O^{p'}(G_{x,y}) \leq G_y^{[1]}$ and so $ O^{p'}(G_{x,y}) =O^{p'}(G_y^{[1]})$ is normal in $G_y$.
Because $q=p\in \{3,7\}$, $O^{p'}(G_x)/Q_x \cong \PSL_3(p)$ or $\SL_3(p)$ and, as the Schur multiplier of $\PSL_2(p)$ has order $2$, we get $O^{p'}(G_y^{[1]})/Q_x \cong  \ASL_2(p)$. In particular, $G_x^{[2]} \cap O^{p'}(G_y^{[1]})=Q_x$.

Let $a \in  O^{p'}(G_y^{[1]})$ be an involution and set $T=C_y\langle a\rangle $.  Then  \cref{lem:P5} (iv)  implies $$[T^{\Delta (x)}, O^{p'}(G_y^{[1]})^{\Delta (x)}]= [a^{\Delta (x)},O^{p'}(G_y^{[1]})^{\Delta (x)}] [C_y^{\Delta (x)}, O^{p'}(G_y^{[1]})^{\Delta (x)}]\leq Q_y^{\Delta(x)}$$ and so   $T\le C_yG_x^{[1]}$, as $|C_y^{\Delta (x)}|$ is even. Hence using \cref{lem:P2} (i) $$ T = C_yG_x^{[1]}\cap  T= C_y( G_x^{[1]}\cap T)\le C_y (G_x^{[1]}\cap G_y^{[1]})=C_yG_x^{[2]}=C_y.$$ Thus $a\in C_y\cap O^{p'}(G_y^{[1]})$. Since $O^{p'}(G_y^{[1]})$ is normal in $G_y$,  we deduce that $S_y \leq O^{p'}(G_y^{[1]})$ which is impossible as $O^{p'}(G_y^{[1]})\cap G_x^{[2]}= Q_x$. This completes the proof.
\end{proof}

\begin{lemma}\label{lem:L1}
Suppose that $G_x^{\Delta(x)}\cong \PSL_3(2)$. Pick $U \in {\mathrm Syl}_3(G_y)$ and set $D= U\cap C_y$ and $F =D \cap L_{x,y}$. Then
\begin{enumerate}[(i)]
\item $q_x=6$ and $q_y=2$;
\item $|G_x^{[2]}/Q_x|=3$;
\item   $D$ is elementary abelian of order $9$,  $C_y= DQ_y$ and  $D \in {\mathrm Syl}_3(G_y^{[1]})$;
\item  $F$ has order $3$, $Q_yF \trianglelefteq G_y$, $F=Z(U)$ and $U$ is extraspecial of order $27$.
\end{enumerate}
\end{lemma}

\begin{proof}
We have $p=2$, $q_x=6$ and \cref{lem:A2.5} implies  $|G_x^{[2]}/Q_x|=3$.

   By \cref{lem:Cy}, $C_y/Q_y$ is an elementary abelian $3$-group of rank at least $2$.
   Since $G_{x}^{[1]} \cap G_{y}^{[1]} \le G_{x}^{[2]}$ and $G_{x,y}^{\Delta(x)} \cong \Sym(4)$, we deduce $C_y/Q_y$ has order $9$. Hence $D$ is elementary abelian of order $9$, $C_y= DQ_y$ and $D \in {\mathrm Syl}_3(G_y^{[1]})$. This proves  (iii).

We know $FQ_y$ has index $3$ in $C_y$ and $FQ_y$ is normalized by $G_{x,y}$. If $FQ_y= G_{u}^{[2]}Q_y$ for some $u \in \Delta(y)$, then as $G_{x,y} $ is transitive on $\Delta(y)\setminus\{x\}$, we have $G_{u}^{[2]}Q_y=G_{z}^{[2]}Q_y$ contrary to \cref{lem:P3} (iv).  So of the four subgroups of index $3$ in $C_y$, there are only three candidates for $G_{u}^{[2]}Q_y$ and so we conclude that $q_y=2$ and this proves (i). Furthermore, we have $FQ_y$ is normalized by $G_y$ and $U$ acts transitively on $\Delta(y)$ and so permutes the three subgroups $G_u^{[2]} Q_y \le C_y$, $u \in \Delta(y)$ transitively. In particular, as $D \in {\mathrm Syl}_3( G_y^{[1]})$ and $G_y/G_y^{[1]}\cong \Sym(3)$, we have $U$ is non-abelian of order $27$ and $F = Z(U)$.
\end{proof}

\begin{lemma}\label{lem:nccf Gx} Assume that $G_x^{\Delta(x)}\cong \PSL_3(2)$. Then
$L_x$ has    either $1$, $2$, $3$ or $6$ non-central $L_x$ chief factors each of which is $3$-dimensional.
\end{lemma}

\begin{proof} We establish \cite[Hypothesis]{MS} with $p=2$ using boldface letters for the groups used in \cite[Hypothesis]{MS}.   So set ${\mathbf M}= L_x$, ${\mathbf E}= Q_x$, ${\mathbf B}\in {\mathrm Syl}_2(\mathbf M)$ and $\mathbf P_1, \mathbf P_2\le \mathbf M$ such that $\mathbf P_1 \cap \mathbf P_2=\mathbf B$.  To reassure ourselves, this means that ${\mathbf M}/{\mathbf E}=L_x/Q_x \cong \PSL_3(2)$, ${\mathbf B}/{\mathbf E} \cong {\mathrm Dih}(8)$ and ${\mathbf P_1}/{\mathbf E}\cong {\mathbf P_2}/{\mathbf E}\cong \Sym(4)$. We choose notation  so that ${\mathbf P_1}=L_{x,y}$. We have $\mathbf P_1=\mathbf P_1^*$ and $\mathbf P_2=\mathbf P_2^*$. This means that \cite[Hypothesis (WBN)]{MS} is satisfied. So we take $\mathbf S=\mathbf B$ and ${\mathbf T}= Q_y$.  By \cref{lem:no char sbgrp}, \cite[Hypothesis (P)]{MS} holds. Since $\mathbf P_1=\mathbf P_1^*$ and $\mathbf P_2=\mathbf P_2^*$, $O^{2'}(\mathbf P_1)= \mathbf P_1$ and $O^{2'}(\mathbf P_2)= \mathbf P_2$. Thus, setting
$$\mathbf L=\langle O^2(O^{2'}(\mathbf P_1^*)),O^2(O^{2'}(\mathbf P_2^*))\rangle= \langle O^2(\mathbf P_1^*),O^2(\mathbf P_2^*) \rangle$$ and remembering $\mathbf E \le \mathbf T$, we have $$\mathbf L \mathbf E= \mathbf L \mathbf T.$$ So, as $\mathbf E$ is a $2$-group, $\mathbf E/  O_2(\mathbf E)$ is a $2'$-group, and $\mathbf E\le \mathbf FO_2(\mathbf E)$ for every subgroup $\mathbf F \le \mathbf M$ with $$\mathbf L  \le \mathbf E\mathbf F.$$ Hence \cite[Hypothesis (A)]{MS} and   \cite[Hypothesis (B)]{MS} both hold.

Since \cite[Hypothesis]{MS} holds, and since ${\mathbf M}= L_x$, we can conclude from \cite[Theorem]{MS} that one of the cases (4), (5), (8), (11) or (13) of that theorem hold. In particular, we see that
${\mathbf E}= Q_x$ has either $1$, $2$, $3$ or $6$ non-central $L_x$ chief factors each of which is $3$-dimensional. This proves the claim.
\end{proof}

\begin{lemma} \label{lem:step2}
\begin{enumerate}[(i)]\item
$L_x$ has     two  non-central   chief factors in $Q_x$ each of which is $3$-dimensional.
\item $F$ has three non-central   chief factors in $Q_y$ (where $F$ is as in \cref{lem:L1}).
\end{enumerate}
\end{lemma}

\begin{proof} Consider the non-central chief factors of $G_y$ in $Q_y$.  If the chief factor is not centralized by $F= Z(U)$, then, as $U$ is extraspecial of order $27$, the chief factor has order a multiple of $2^6$ and $F$ acts fixed-point-freely.  Thus the number of $F$-chief factors is a multiple of $3$ and  which we denote by $3f$. As $G_y$ has characteristic $p$, $f \ge 1$.

From the perspective of $L_x$, we have $Q_y/Q_x$ is a non-central $F$-chief factor,  and, as each $L_x$ non-central chief factor in $Q_x$ is $3$-dimensional by  \cref{lem:nccf Gx},  $FQ_x$ has one non-central chief factor for each $L_x$ non-central chief factor. Thus $L_x$ has $3f-1$ non-central chief factors in $Q_x$. Using \cref{lem:nccf Gx}, we deduce that $L_x$ has $2$ non-central chief factors in $Q_x$. Thus (i) holds and (ii) follows from this.
\end{proof}

Recall the definition of $Z_x$ and $Z_y$ from \cref{not:1stepmore} and
observe $Z_y \leq Z_x$, since $\Delta$ is of local characteristic $2$ and $Q_x \leq Q_y$.

\begin{lemma}\label{lem:Z1} Suppose $G_x^{\Delta(x)}\cong \PSL_3(2)$. Then
\begin{enumerate} [(i)]
\item $Q_y=Q_xQ_z$ and  $Z_x \cap Z_z=Z_y$;
\item $[Z_x, F] \neq 1$.
\end{enumerate}
\end{lemma}

\begin{proof}
We continue the notation from \cref{lem:L1}.

$(i)$. We know $FQ_y \le G_y^{[1]}$. Since $F$ normalizes $Q_z$ and $Q_z \ne Q_x$, we have $Q_xQ_z= Q_y$.
Since $\Delta$ is of local characteristic $p$ and $Q_y=Q_xQ_z$,
we glean $Z(Q_y) \leq Z(Q_x) \cap Z(Q_z) \leq Z(Q_y)$.
\ms
\par\noindent
$(ii)$.
Suppose $[Z_x, F]=1$. Then $[Z_x, O^2(L_x)] = 1$  and
$L_x=Q_xO^2(L_x)$. Since $O^2(L_x)$ acts transitively on $\Delta(x)$ and
$Z_y\leq Z_x$, we obtain $Z_y=1$ from \cref{lem:basic1} and this is a contradiction.\end{proof}

\begin{lemma}\label{lem:zznotin} Suppose $G_x^{\Delta(x)}\cong \PSL_3(2)$. Then $Z_z\not \le G_x^{[1]}$.
\end{lemma}
\begin{proof} Again we use the notation started in \cref{lem:L1}. By \cref{lem:A2.5}, we know $G_x^{\Delta(x)}\cong \PSL_3(2)$.
Set $V_y = \langle Z_u\mid u \in \Delta(y)\rangle$. Then $V_y\leq Q_y$ is normalized by $G_y$ and  $F$ does not centralize $Z_x\le V_y$ by \cref{lem:Z1} (ii). Hence $U$ acts faithfully on $V_y$ and so we conclude that $F$ has at least three non-central chief factors in $V_y$. Since $F$ has exactly three non-central chief factors in $Q_y$ and $Q_y/Q_x$ is such a factor, we conclude that $V_y \not \le G_x^{[1]}$ and this delivers the claim.
\end{proof}

\begin{proof}[Proof of \cref{prop:notlnq}] \cref{lem:A2.5} gives $G_x^{\Delta(x)}\cong \PSL_3(2)$ and, by \cref{lem:zznotin}, $Z_z \not\leq Q_x$.  As $F\leq G_y^{[1]}$ normalizes $Z_z$, $Q_y=Q_xZ_z=Z_xQ_z$
Since $Q_x\cap Z_z $ is centralized by $Z_xQ_z$, we deduce $Q_x \cap Z_z=Z_y$.

Let $u \in \Delta (y) \setminus \{x,z\}$. Since $s \geq 4$, we have $Z_z \not\leq Q_{u}$. We calculate using the fact that $Q_x$ is the unique Sylow $2$-subgroup of $G_{x}^{[2]}$ by \cref{lem:L1}  (ii)$$[G_{x}^{[2]},Z_z]\leq G_x^{[2]} \cap Z_z \leq Q_x \cap Z_z=Z_y$$ and
$$[G_u^{[2]},Z_z]\leq Q_{z^\prime} \cap Z_z=Z_y.$$
It follows that
$$[F ,Z_z] \le [G_u^{[2]}G_x^{[2]},Z_z] \leq Z_y\le Q_x.$$ However, $Z_zQ_x/Q_x$ is not centralized by $F$, and so this is impossible. This contradiction completes the proof.
\end{proof}

\section{The main theorem}

Suppose that $X$ is a $2$-transitive group in its action on $\Omega$.  Then \cite[Lemma 2.2]{bonstell} (for example) yields that either there is a prime $r$ such that $F^*(X)$ is a regular elementary abelian $r$-group, or $F^*(X)$ is a non-abelian simple group. In the first case, we say that $X$ is of \emph{regular type} and in the second  that $X$ is of \emph{simple type}. When $X$ is of simple type, the description of $F^*(X)$ and $\Omega$ is conveniently presented in \cite[Lemma 2.5]{bonstell} (this result requires the classification of the finite simple groups). Since $G_x^{\Delta(x)}$ acts $2$-transitively on $\Delta(x)$ and as we also know that $1\ne Q_y^\Delta(x) \unlhd G_{x,y}^{[1]}$, \cref{prop:notlnq} combined with \cite[Lemma 2.5]{bonstell} yields

\begin{lemma}\label{lem:A1}  The group
$F^*(G_x^{\Delta (x)})$ is either of regular type or is of simple type and is isomorphic to a rank 1 group of Lie type
in characteristic $p$
in its natural permutation representation (including $\mathrm{Ree}(3)^\prime$).
\end{lemma}\qed

 We can now move directly to the proof \cref{thm:main}.

\begin{proof}[Proof of \cref{thm:main}] Set $N= O^p(L_x)Q_x$. Then $O_p(N)= Q_x$ and $N=  O^p(N)Q_x$. Define $S= Q_y$ and $\widetilde G= G_x/Q_x$. Then with $L= L_x$, we have $L= NS$, $O_p(N)= O_p(L)=Q_x < Q_y=S$ and $C_S(Q_x) \le C_G(Q_x)\le Q_x$. Furthermore, \cite[Hypothesis 3.3(b) and (c)]{bonstell} follows from \cref{lem:P1}(v) and \cref{lem:A1} respectively. Because of \cref{lem:no char sbgrp} and \cite[Lemma 3.8]{bonstell}  we are in the same conclusions as  \cite[7.7 steps $\mathbf 1^\circ$ and $\mathbf 2^\circ$]{bonstell}. Following \cite[7.7 steps $\mathbf 3^\circ$ through $\mathbf 9 ^\circ$]{bonstell} verbatim (being careful to note the role of $x$ and $y$ are reversed) yields
$O^p(L_x) \cong \AGL_2(q)'$, $q=p^r$ with $p$ odd, and $Q_x$ elementary abelian.   This completes the proof.
 \end{proof}

\bigskip\noindent
{\bf Statements and Declarations}
\par\smallskip
\begin{itemize}
\item The authors have no conflicts of interest to declare.
\item Data sharing not applicable to this article as no datasets were generated or analysed during the current study.
\end{itemize}

\vspace{1cm}

\noindent
{{\em Mathematics Subject Classification (2010)}  20B25, 05C25, 05E18\\ Keywords: locally $s$-arc transitive graphs, group amalgams}

\bigskip\noindent
\footnotesize\noindent
\textsc{Dipartimento di Matematica e Informatica, Universit\`a della Calabria,
    87036 Arcavacata di Rende, Italy}\par\nopagebreak\noindent
  \textit{E-mail address}: \texttt{vanbon@mat.unical.it}

\smallskip\noindent
\textsc{School of Mathematics, University of Birmingham, Edgbaston,
    Birmingham B15 2TT, UK}\par\nopagebreak\noindent
  \textit{E-mail address}: \texttt{c.w.parker@bham.ac.uk}

\end{document}